\documentclass{elsarticle}
\usepackage{amsmath}
\usepackage{hyperref}
\usepackage{cleveref}
\usepackage{color}
\usepackage{graphicx}
\graphicspath{{IMAGES/}}
\usepackage[bb=px, bbscaled=1,cal=euler,scr=boondoxo,scrscaled=1.04,frak=euler,frakscaled=1.05]{mathalpha}

\usepackage{cprotect}  

\mathcode`\@="8000
{\catcode`\@=\active \gdef@{\mkern1mu}}

\def\mydate{\number\day\ {\ifcase\month \or January\or February\or
              March\or April\or May\or June\or July\or August\or
              September\or October\or November\or December\fi}
\number\year}

\def\whiteout#1#2{\color{white}{\rule{#1}{#2}}}

\usepackage{colonequals}
\usepackage{mathdots}

\def\Lyap{{\boldsymbol{\cal L}}}
\def\Id{{\boldsymbol{\cal I}}}

\def\BA{{\bf A}}  
\def\Bb{{\bf b}}\def\BB{{\bf B}}  
\def\BC{{\bf C}}  
\def\BD{{\bf D}}

\def\BG{{\bf G}}  
  
\def\BI{{\bf I}}  
  
\def\BK{{\bf K}}  \def\CK{{\cal K}}
  
\def\BM{{\bf M}}

  \def\CP{{\cal P}}
  
\def\Br{{\bf r}}\def\BR{{\bf R}}

\def\Bv{{\bf v}}\def\BV{{\bf V}}  
\def\Bw{{\bf w}}\def\BW{{\bf W}}  
\def\Bx{{\bf x}}\def\BX{{\bf X}}

\def\BAs{\BA\kern-1.5pt}
\def\Bzero{\boldsymbol{0}}
\def\BLambda{\boldsymbol{\Lambda}}

\def\R{\mathbb{R}}
\def\C{\mathbb{C}}

\def\Cn{\C^n}

\def\Cnn{\C^{n\times n}}

\def\eop{{\rm e}}

\def\Re{{\rm Re}}
\def\hermA{{\textstyle{\frac{1}{2}}(\BA+\BAs^*)}}

\def\ipG#1{\langle #1\rangle_\BG}
\def\normG#1{\|#1\|_\BG}

\def\CPs{\CP{\kern-0.8pt}}
\def\mingmres{\min_{\stackrel{\scriptstyle{p\in \CPs_k}}
                            {\scriptstyle{p(0)=1}}}}

\newtheorem{thm}{Theorem}
\newtheorem{cor}{Corollary}

\begin{document}

\title{Extending Elman's Bound for GMRES}
\author[1]{Mark Embree}
\ead{embree@vt.edu}
\affiliation[1]{organization={Department of Mathematics, Virginia Tech},
                addressline={\\[2pt] 225 Stanger Street 0123},
                city = {Blacksburg, Virginia},
                postcode = {24061},
                country = {USA}}

\begin{abstract}
If the numerical range of a matrix is contained in the right half of the
complex plane, the GMRES algorithm for solving linear systems will reduce
the norm of the residual at every iteration.  
In his Ph.D.\ dissertation, Howard Elman derived a bound that
guarantees convergence.  
When the numerical range contains the origin, GMRES need not make progress 
at every step and Elman's bound does not apply, even if all the eigenvalues
are located in the right half-plane.  
However by solving a Lyapunov equation, one can construct an inner product 
in which the numerical range is contained in the right half-plane.
One can then bound GMRES (run in the standard Euclidean norm) by applying
Elman's bound in this new inner product, at the cost of a multiplicative 
constant that characterizes the distortion caused by the change of inner 
product.  
Using Lyapunov inverse iteration, one can build a family of 
such inner products, trading off the location of the numerical
range with the size of constant. 
This approach complements techniques that
Greenbaum and colleagues have recently proposed for excising the origin 
from the numerical range to gain greater insight into the convergence
of GMRES for nonnormal matrices.
\end{abstract}


\maketitle


\section{Introduction}

Suppose the eigenvalues of $\BA\in\Cnn$ are all located in the open right half-plane.
We seek to solve the system $\BA\Bx=\Bb$ using the GMRES algorithm~\cite{SS86},
which produces at its $k$th iteration a solution estimate $\Bx_k$ that minimizes
the 2-norm of the residual $\Br_k := \Bb-\BA\Bx_k$ over all vectors in the 
$k$th Krylov subspace $\CK_k(\BA,\Bb) = {\rm span}\{\Bb,\BA\Bb,\ldots,\BAs^{k-1}\Bb\}$.
Various strategies exist for bounding $\|\Br_k\|_2$ as a function of $k$, based on
spectral properties of $\BA$; see \cite{Emb22} for a survey.
A bound proposed by Howard Elman~\cite[thm.~5.4]{Elm82} (see also \cite{EES83})
has particular appeal because of its simplicity.  
Let $\mu(\BA)$ denote the leftmost eigenvalue of 
$\frac{1}{2}(\BA+\BAs^*)$, the Hermitian part of $\BA$.\ \ 
Provided $\mu(\BA)>0$, Elman showed that
\begin{equation} \label{eq:elman}
 \frac{\|\Br_k\|_2}{\|\Bb\|_2} \le  \bigg(1 - \frac{\mu(\BA)^2}{\|\BA\|_2^2}\bigg)^{\!k/2}.
\end{equation}
Beckermann, Goreinov, and Tyrtyshnikov~\cite{BGT06} explained 
how this bound could be understood in terms of the numerical range
(or field of values) of $\BA$,
\begin{equation} \label{eq:W}
 W(\BA) = \left\{\frac{\Bv^*\BA\Bv}{\Bv^*\Bv}: \Bzero \ne \Bv\in\Cn\right\}.
\end{equation}
The numerical range
is a closed, convex subset of the complex plane that contains
the eigenvalues of $\BA$.\ \    
(For these and other properties of $W(\BA)$, see~\cite[chap.~1]{HJ91}.)
For $z := \Bv^*\BA\Bv/\Bv^*\Bv \in W(\BA)$, 
\[ \Re(z) = {\textstyle{\frac{1}{2}}} (z + \overline{z}) 
          = \frac{1}{2} \bigg(\frac{\Bv^*\BA\Bv}{\Bv^*\Bv} 
               + \frac{\Bv^*\BAs^*\Bv}{\Bv^*\Bv}\bigg)
          = \frac{\Bv^* \big(\hermA\big)\Bv}{\Bv^*\Bv}.\]
Thus the real part of any point $z$ in the numerical range is a Rayleigh quotient 
for the Hermitian part of $\BA$.  
Taking extremes of this Rayleigh quotient, we have
\[ \Re(W(\BA)) = \Big[ \lambda_{\min}(\hermA), \lambda_{\max}(\hermA)\Big].\]
Thus, $\mu(\BA)>0$ is necessary and sufficient for $W(\BA)$ to be contained
in the right half-plane.
Any $z = \Bv^*\BA\Bv/\Bv^*\Bv \in W(\BA)$ satisfies $|z| \le \|\BA\|_2$, so
the numerical range can be bounded in a disk segment determined by $\mu(\BA)$ and $\|\BA\|_2$:
\begin{equation} \label{eq:nrbound}
 W(\BA) \subseteq \{z \in \C: \Re(z)\ge \mu(\BA)\}\ \cap\ \{z\in\C: |z|\le \|\BA\|_2\}.
\end{equation}

\begin{figure}[b!]
\begin{center}
\includegraphics[height=1.8in]{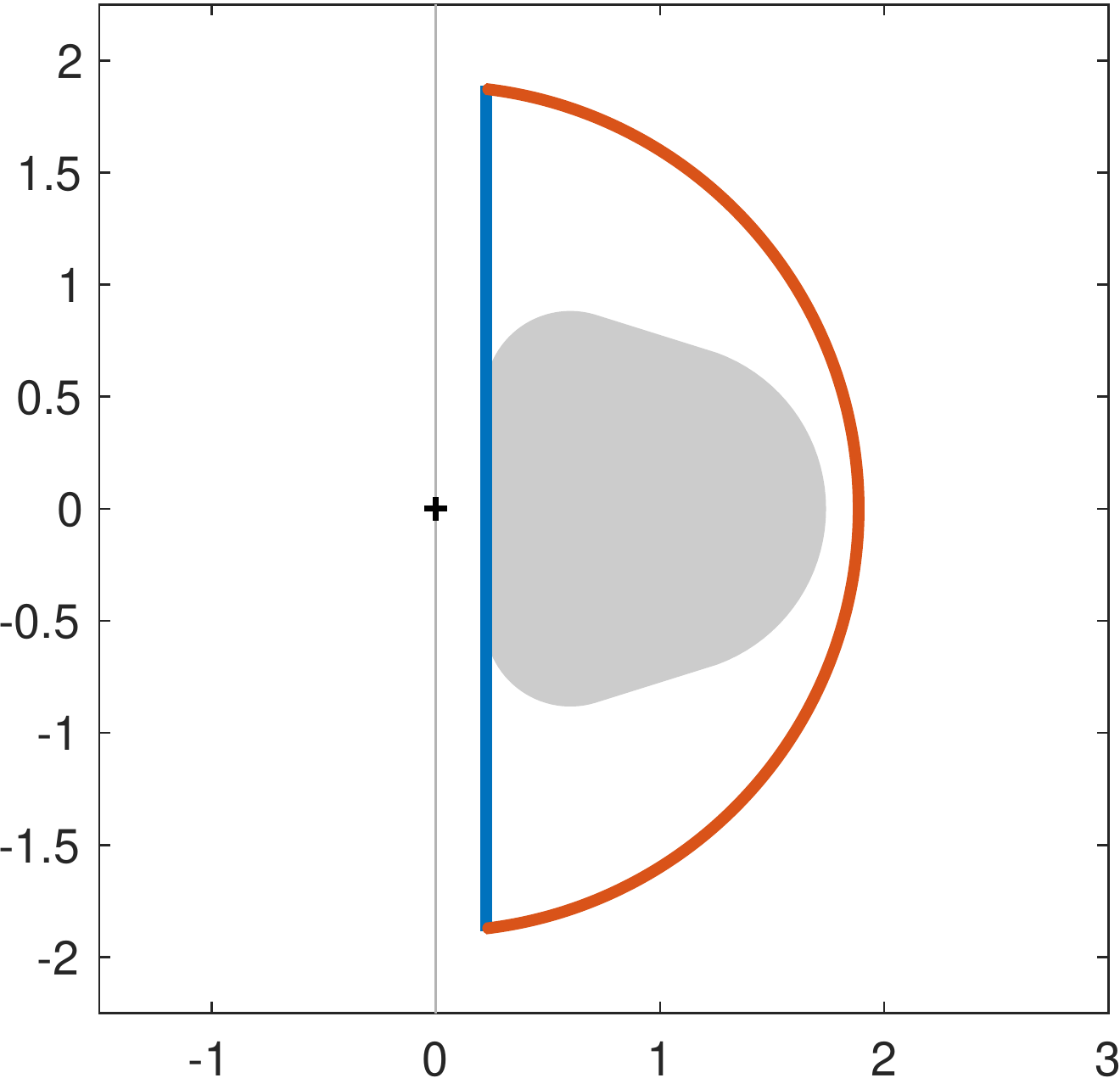}
\begin{picture}(0,0)
\put(-69,65){\footnotesize $W(\BA)$}
\end{picture}
\end{center}

\vspace*{-8pt}
\caption{\label{fig:cartoon}
Illustration of the bound~\cref{eq:nrbound}, in which 
the numerical range $W(\BA)$ is contained in a disk segment: 
the left edge of the disk segment (blue) is determined by $\mu(\BA)$, 
the leftmost eigenvalue of the Hermitian part of $\BA$; 
this bound is sharp, i.e., there exists $z\in W(\BA)$ with ${\rm Re}(z)=\mu(\BA)$. 
The circular arc (red) is determined by $\|\BA\|_2$; it need not be sharp, 
but cannot be off by more than a factor of two: there exists $z\in W(\BA)$
such that $|z|\ge \frac{1}{2}\|\BA\|_2$~{\rm \cite[(5.7.21)]{HJ13}}.}
\vspace*{-4pt}
\end{figure}

\Cref{fig:cartoon} illustrates this set.  
Beckermann et al.~\cite{BGT06} observed that their analysis applies 
to any matrix $\BA\in\Cnn$ whose numerical range does not include the origin, 
as one can rotate the problem to $\eop^{\iop@\theta}\BA \Bx = \eop^{\iop@\theta}\Bb$
for some $\theta\in[0,2\pi)$ so that $\mu(\eop^{\iop@\theta}\BA)>0$.  
The same observation applies to the analysis that follows.

To study transient growth of solutions $\Bx(t) = \eop^{t\BA}\Bx(0)$ 
to the evolution problem $\Bx'(t)=\BA\Bx(t)$ 
(and its discrete time analog $\Bx_{k+1} = \BA\Bx_k$),
Godunov~\cite[sect.~9.3]{God97} proposed a bound on $\|\eop^{t\BA}\|$ 
that, in essence, finds an inner product in which 
the numerical range of $\BA$ is contained
in the left half-plane (thus ensuring monotone convergence 
of the induced norm of $\eop^{t\BA}$ as $t\to\infty$)~\cite[p.~147]{TE05plain}; 
for various explications and applications of this bound, 
see~\cite{AAM23,Gar23,Pli05,Ves03a,Ves11}.  

This note applies Godunov's idea to the GMRES setting, to develop a framework
in which to think about (and generate) alternative inner products for 
analyzing GMRES convergence.  This work complements the recent analysis
of Spillane and Szyld for ``weighting'' (changing inner products), 
``deflating'', and ``preconditioning'' GMRES~\cite{SS24}.
Here we observe that, for systems where the spectrum is contained on one side
of the origin,  one can solve a Lyapunov equation to get an 
inner product for which the numerical range of $\BA$ does not contain 
the origin, and apply Elman's estimate (or the Crouzeix--Palencia 
theorem~\cite{CP17}, or indeed various other GMRES bounds) 
in this inner product.
The price for applying a GMRES bound in this new inner product
is a multiplicative constant that describes the conditioning of 
the Gram matrix for the inner product, 
that is, how much the new inner product distorts the original geometry.

We mainly view this approach as a theoretical tool for studying GMRES convergence for
nonnormal matrices,
although we describe several applications where the Lyapunov equation
has an explicit solution (including a preconditioned saddle point problem
that is positive definite in a ``nonstandard inner product''~\cite{FRSW98,SW08}).
In addition, this general framework could potentially inform the design of nonnormality-reducing 
preconditioners for special classes of problems~\cite{GL01} or
GMRES variants that adaptively change inner products upon each 
restart~\cite{EMN17,Ess98,GP14}.

\section{Finding a favorable inner product}

Let $\BC\in\Cnn$ be \emph{any} Hermitian positive definite matrix.
Since all eigenvalues of $\BA$ are in the right half-plane, 
the Lyapunov matrix equation
\begin{equation}\label{eq:lyap}
    \BAs^*\BG + \BG \BA = \BC
\end{equation}
has a unique solution $\BG\in\Cnn$, and it is Hermitian positive definite;
see, e.g., \cite[thm.~2.2.3]{HJ91}.
Take its Cholesky (or square root) factorization, $\BG = \BR^*\BR$.\ \ 
Use $\BG$ to define the inner product
\[ \ipG{\Bv,\Bw} := \Bw^*\BG \Bv = (\BR\Bw)^*(\BR\Bv)\]
with the associated norm
\[ \normG{\Bv} := \sqrt{\ipG{\Bv,\Bv}} = \|\BR\Bv\|_2.\]
The matrix norm induced by this vector norm satisfies,
for any $\BX\in\Cnn$,
\[ \normG{\BX} 
   := \max_{\Bv\ne \Bzero} \frac{\normG{\BX\Bv}}{\normG{\Bv}} 
    = \max_{\Bv\ne \Bzero} \frac{\|\BR\BX\BR^{-1}\BR\Bv\|_2}{\|\BR\Bv\|_2} 
    = \max_{\Bw\ne \Bzero} \frac{\|\BR\BX\BR^{-1}\Bw\|_2}{\|\Bw\|_2}
    = \|\BR\BX\BR^{-1}\|_2.
\]

\smallskip
\noindent
The matrix $\BR$ can be computed in a few lines of MATLAB code.\\[4pt]
\hspace*{15pt} \verb|L = chol(C)';          % Cholesky factorization of C|\\
\hspace*{15pt} \verb|R = lyapchol(-A',L);   % Lyapunov solution in factored form| \\
\hspace*{15pt} \verb|A_G = R*(A/R);         % Note: norm(A_G) = G-norm of A| \\[4pt]
The Lyapunov solver \verb|lyapchol| calls the SLICOT~\cite{BMSVV99} routine \verb|SB03OD| 
to compute the upper-triangular factor $\BR$ for the solution $\BG=\BR^*\BR$ directly via
Hammarling's algorithm~\cite{Ham82}, ensuring the positive definiteness 
of the computed solution.

\medskip
We claim that the numerical range of $\BA$ \emph{in the $\BG$-inner product},
\[ W_\BG(\BA) := \left\{\frac{\ipG{\BA\Bv,\Bv}}{\ipG{\Bv,\Bv}}: 
                   \Bzero\ne \Bv\in\Cn\right\}
               = \left\{\frac{\Bv^*\BG\BA\Bv}{\Bv^*\BG \Bv}: 
                   \Bzero\ne \Bv\in\Cn\right\},
\]
is contained in the right half-plane.
To see this, take any $z =  \ipG{\BA\Bv,\Bv}/\ipG{\Bv,\Bv} \in W_\BG(\BA)$,
and note that
\begin{eqnarray*} 
\Re(z) = {\textstyle{\frac{1}{2}}} (\overline{z}+z) 
          &=& \frac{1}{2} \bigg(\frac{\Bv^*\BAs^*\BG\Bv}{\Bv^*\BG\Bv} 
                           + \frac{\Bv^*\BG\BA\Bv}{\Bv^*\BG\Bv}\bigg) \\[3pt]
          &=& \frac{\Bv^* \big(\BAs^*\BG+\BG\BA\big)\Bv}{2@\Bv^*\BG\Bv} 
           =  \frac{\Bv^* \BC\Bv}{2@\Bv^*\BG\Bv} > 0,
\end{eqnarray*}
since $\BG$ solves the Lyapunov equation~\cref{eq:lyap},
and both $\BG$ and $\BC$ are positive definite.

\section{GMRES bounds}
With the $\BG$-norm in hand, we are prepared to bound the 2-norm of the
GMRES residual at the $k$th iteration.  Recall that this residual satisfies 
\begin{equation} \label{eq:bound0}
 \|\Br_k\|_2 = \mingmres \|p(\BA)\Bb\|_2,
\end{equation}
where $\CP_k$ denotes the set of polynomials of degree $k$ or less.
Then
\begin{eqnarray*}
   \|\Br_k\|_2 &=& \mingmres \|\BR^{-1} \BR@ p(\BA) \BR^{-1}\BR\Bb\|_2 \\[3pt]
 &\le& \mingmres \|\BR^{-1}\|_2 \|\BR @p(\BA) \BR^{-1}\|_2 \|\BR\|_2 \|\Bb\|_2 \\[3pt]
 &=& \mingmres \|\BR^{-1}\|_2 \|\BR\|_2\,\normG{@p(\BA)} \|\Bb\|_2.
\end{eqnarray*}
Noting that $\|\BG\|_2 = \|\BR^*\BR\|_2 = \|\BR\|_2^2$ and 
$\|\BG^{-1}\|_2 = \|\BR^{-1}\BR^{-*}\|_2 = \|\BR^{-1}\|_2^2$,
we can rearrange the bound as follows.

\begin{thm} \label{thm:main}
Let $\BA\in\Cnn$ be a matrix having all eigenvalues in the open right half-plane.
Given any Hermitian positive definite $\BC\in\Cnn$, let $\BG$ denote the solution
of the Lyapunov equation $\BAs^*\BG + \BG\BA = \BC$.  
At the $k$th iteration, the residual $\Br_k$ produced by GMRES 
applied to $\BA\Bx=\Bb$ with starting vector $\Bx_0=\Bzero$ satisfies
\begin{equation} \label{eq:main}
   \frac{\|\Br_k\|_2}{\|\Bb\|_2} \,\le\, \sqrt{\kappa_2(\BG)}\,\mingmres \normG{@p(\BA)},
\end{equation}
where $\kappa_2(\BG) := \|\BG\|_2 \|\BG^{-1}\|_2$.
\end{thm}

The bound~\cref{eq:main} has the form and spirit of bounds
recently proposed by Chen, Greenbaum, and Wellen~\cite{CGW}, 
in which a similarity transformation $\BR\BA\BR^{-1}$
is applied to $\BA$ to yield a matrix more amenable to analysis, 
at the cost of the multiplicative constant
$\kappa_2(\BR)$. 
The contribution of Theorem~\ref{thm:main} is use of the Lyapunov equation as
a means for identifying such similarity transformations, for this class 
of matrices $\BA$.\ \ 
(See~\cite[chap.~51]{TE05plain} for a discussion of how similarity 
transformations affect the pseudospectra of $\BA$.)

\medskip
Since $W_\BG(\BA)$ is contained in the open right half-plane, a variety
of approaches can be applied to~\cref{eq:main} to obtain more explicit
convergence bounds. 
Elman's theorem~\cref{eq:elman}
only uses data that bounds the numerical range in the intersection of a 
half-plane and a disk, as seen in \cref{fig:cartoon}.
An analogous containment holds in the $\BG$-inner product: 
\[ W_\BG(\BA) \subseteq 
   \{z\in\C: \Re(z) \ge \mu_\BG(\BA)\}
   \ \cap\ 
   \{ z \in \C: |z| \le \|\BA\|_\BG\},\]
where $\mu_\BG(\BA) = \min_{z\in W_\BG(\BA)} \Re(z) 
                    = \lambda_{\min}\big({\textstyle{\frac{1}{2}}}(\BR\BA\BR^{-1}+\BR^{-*}\BAs^*\BR^*)\big)$,
and $\|\BA\|_\BG = \|\BR\BA\BR^{-1}\|_2$.
We can apply Elman's bound in the $\BG$-inner product to~\cref{eq:main} in \Cref{thm:main}.

\begin{cor} \label{cor:elman}
Under the conditions of \Cref{thm:main},
\begin{equation} \label{eq:new_elman}
   \frac{\|\Br_k\|_2}{\|\Bb\|_2} 
     \,\le \, \sqrt{\kappa_2(\BG)}\ \bigg(1 - \frac{\mu_\BG(\BA)^2}{\|\BA\|_\BG^2}\bigg)^{\!k/2}.
\end{equation}
\end{cor}

Using the same ratio $\mu(\BA)/\|\BA\|_2$,
Beckermann, Goreinov, and Tyrtyshnikov~\cite[thm.~2.1]{BGT06} 
derived a bound with an improved asymptotic convergence rate, 
at the cost of an additional multiplicative constant.
In~\cite[cor.~1.3]{Bec05} Beckermann simplified this constant, 
which is the version given below, cast here in the $\BG$-inner product.%
\footnote{In fact, Beckermann also sharpened the \emph{rate} in this bound 
by replacing $\|\BA\|$ by the numerical radius $\nu(\BA) := \max_{z\in W(\BA)} |z|$, 
which effectively replaces the red arc in~\Cref{fig:cartoon} by a smaller arc that
touches the largest magnitude point in $W(\BA)$.  
Computing $\nu(\BA)$ involves a more sophisticated algorithm~\cite{MO05} than $\|\BA\|$. 
To facilitate comparison with Elman's bound, we use $\|\BA\|$ in~\Cref{cor:bgt}.
We thank Pierre Marchand for pointing out the improvement in~\cite{Bec05}.}

\begin{cor} \label{cor:bgt}
Under the conditions of \Cref{thm:main},
\begin{equation} \label{eq:new_bgt}
   \frac{\|\Br_k\|_2}{\|\Bb\|_2} 
     \,\le\,
          \sqrt{\kappa_2(\BG)}\;(2+\rho_\beta) \,\rho_\beta^k,
\end{equation}
where $\cos(\beta) = \mu_\BG(\BA)/\|\BA\|_\BG$ and $\rho_\beta = 2@\sin\big(\beta\pi/(4\pi-2\beta)\big)$.
\end{cor}

These last two corollaries are based on bounding $W_\BG(\BA)$
by a disk segment, as in \Cref{fig:cartoon}.
Bounding GMRES using $W_\BG(\BA)$ itself can yield a substantially better 
convergence rate.
The Crouzeix--Palencia theorem~\cite{CP17} holds in any inner product, so
\[ \normG{p(\BA)} \le \big(1+\sqrt{2}\big) \max_{z\in W_\BG(\BA)} |p(z)| \]
for all polynomials $p\in\CP_k$, leading to the following bound.

\begin{cor} \label{cor:cp}
Under the conditions of \Cref{thm:main},
\begin{equation} \label{eq:cp}
   \frac{\|\Br_k\|_2}{\|\Bb\|_2} 
     \,\le \,\sqrt{\kappa_2(\BG)}\;(1+\sqrt{2})\,\mingmres \max_{z\in W_\BG(\BA)} |p(z)|.
\end{equation}
\end{cor}

GMRES bounds that use other quantities (e.g., pseudospectra, spectral projectors) 
could also benefit from~\Cref{thm:main}, if the
$\BG$-inner product reduces the departure from normality 
in a way that significantly improves the bounds for $\|p(\BA)\|_\BG$ 
over those for $\|p(\BA)\|_2$.

\section{Examples from applications}

In some applications, a right hand-side $\BC$ and matching solution $\BG$ for 
the Lyapunov equation~\cref{eq:lyap} arise naturally from the structure 
of the problem.  We describe two such examples.
While these settings are quite clean, they suggest forms of $\BG$ that could
be useful in more general contexts, and provide a baseline from which nearby 
problems that lack this exact structure could be analyzed via perturbative
techniques~\cite{SEM13}.

\subsection{A damped mechanical system}

\begin{figure}[t!]
\begin{center}
\includegraphics[height=2.2in]{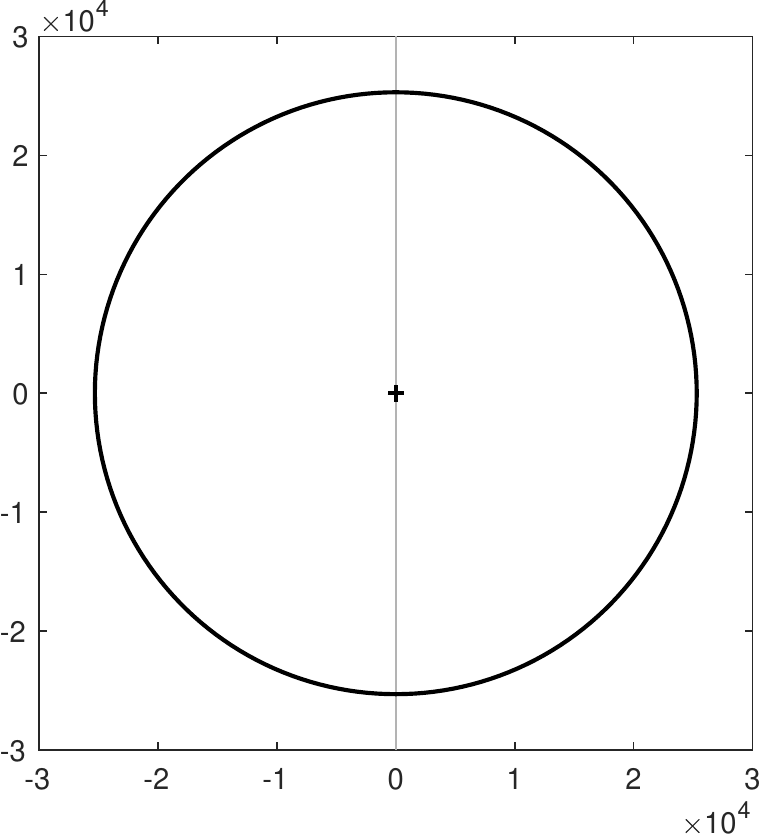}\quad
\raisebox{8pt}{\includegraphics[height=2.02in]{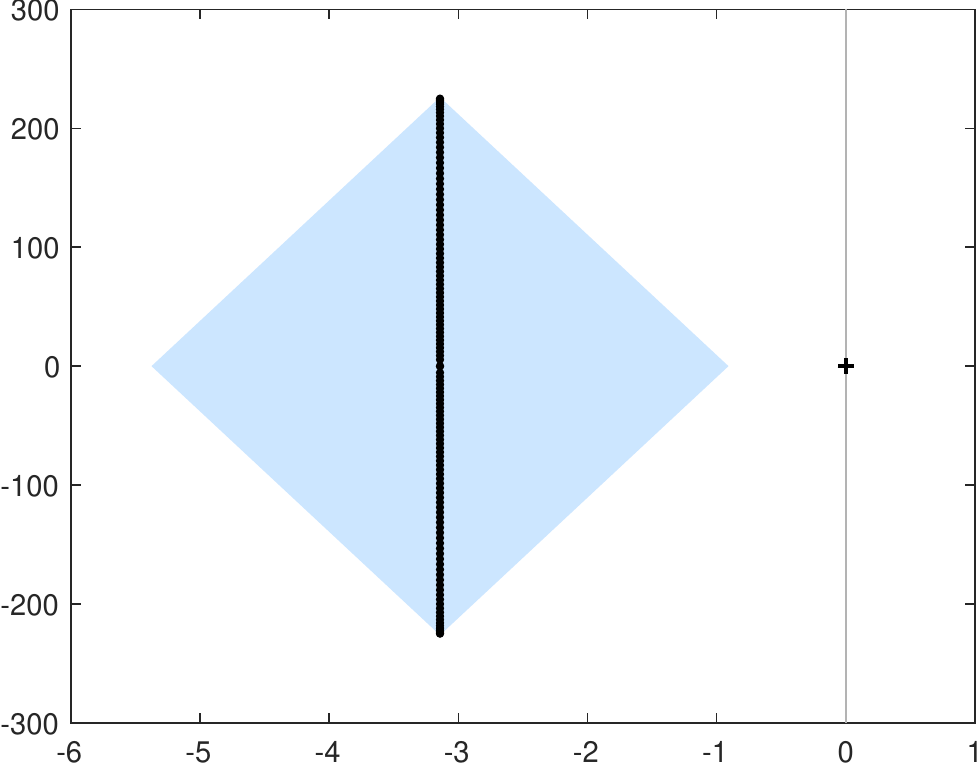}}
\end{center}

\begin{picture}(0,0)
\put(100,103){\footnotesize $W(\BA)$}
\put(255,103){\footnotesize $W_\BG(\BA)$}
\end{picture}

\vspace*{-25pt}
\caption{\label{fig:cox} 
An optimally damped mechanical system~\cref{eq:waveA} for which the Lyapunov
equation~\cref{eq:lyap_cox} has a tidy solution, adapted from~\cite{Cox98a}.
For this example ($N=64$, $n=128$), the origin $(+)$ is embedded deep in the numerical range 
$W(\BA)$.
The inner product defined by $\BG$ in~\cref{eq:waveG} gives a much smaller numerical
range $W_\BG(\BA)$ contained in the left half-plane, with $\smash{\sqrt{\kappa_2(\BG)}} = 225.035$.  
The eigenvalues of $\BA$ (which appear as the black vertical 
segment in the right plot) include a $2\times 2$ Jordan block.  
(The horizontal and vertical axes are scaled differently in the right plot.)}
\end{figure}

Consider the damped mechanical system 
\[\BM \Bx''(t) = -\BK\Bx(t) - \BD \Bx'(t),\]
for state vector $\Bx(t)\in \R^N$, 
Hermitian positive definite mass, stiffness, and damping matrices 
$\BM, \BK, \BD \in\R^{N\times N}$, 
and initial condition $\Bx(0)\in\R^N$.\ \ 
A common problem is the design of $\BD$ to optimally damp the system,
i.e., to drive $\Bx(t)\to\Bzero$ as $t\to\infty$ as fast as possible.
Cox~\cite{Cox98a} considers damping of the form $\BD = 2@a@\BM$,
where $a>0$ is a constant to be optimized.
Writing the second-order differential equation in first-order form leads
to the ``linearization'' matrix
\begin{equation} \label{eq:waveA}
 \BA = \left[\begin{array}{cc} 
             \Bzero & \BI \\[3pt]
             -\BM^{-1}\BK & -\BM^{-1}\BD 
        \end{array}\right]
         = \left[\begin{array}{cc} 
             \Bzero & \BI \\[3pt]
             -\BM^{-1}\BK & -2@a@\BI 
        \end{array}\right].
\end{equation}
In the Euclidean inner product, $W(\BA)$ must contain the origin, 
since the Rayleigh quotient of $\BA$ with $[\Bv^T\ \Bzero^T]^T$ is zero for 
any $\Bv\in\C^N$;
indeed, as seen in~\Cref{fig:cox}, the origin can be embedded far in the 
interior of  $W(\BA)$.\ \ 
For a concrete example, discretize a homogeneous string of unit length, 
giving positive definite
\[ \BM = {\rm tridiag}(h/6,2h/3,h/6), 
    \qquad
   \BK = {\rm tridiag}(-1/h,2/h,-1/h),
\] 
with $h = 1/(N+1)$ for discretization size $N$.
Taking $a$ to be the square root of the 
smallest eigenvalue of $\BM^{-1}\BK$ 
minimizes the spectral abscissa of $\BA$
(giving the fastest asymptotic damping)~\cite[prop.~2.1]{Cox98a},
but this optimal $a$ makes $\BA$ is nondiagonalizable
($-a$ is a double eigenvalue with a $2\times 2$ Jordan block).

To study optimal energy damping, Cox~\cite{Cox98a} notes that the Lyapunov equation
\begin{equation} \label{eq:lyap_cox}
 \BA^*\BG+\BG\BA = -\BC = -\left[\begin{array}{cc} \BK & \Bzero \\ \Bzero & \BM \end{array}\right] 
\end{equation}
has the explicit solution 
\begin{equation} \label{eq:waveG}
\BG = \left[\begin{array}{cc} \frac{1}{2} \BD + \BK\BD^{-1}\BM  & \frac{1}{2}\BM \\[3pt]
         \frac{1}{2} \BM & \BM\BD^{-1}\BM \end{array}\right],
\end{equation}
provided $\BM$, $\BK$, and $\BD$ are (as in this case) simultaneously diagonalizable.
(More precisely, $\BG$ solves the Lyapunov equation provided $\BK\BD^{-1}\BM$ is Hermitian:
$\BK\BD^{-1}\BM = (\BK\BD^{-1}\BM)^* = \BM\BD^{-1}\BK$.)
Since the spectrum of damped wave operators typically falls in the left half-plane,
we take the right-hand side $-\BC$ of~\cref{eq:lyap_cox} to be negative definite;
hence $W_\BG(\BA)$ will be contained in the left half-plane.
\Cref{fig:cox} compares $W(\BA)$ and $W_\BG(\BA)$ for $N=64$ with
$a$ selected to give optimal asymptotic damping -- and hence a defective real eigenvalue.
(Since this $\BA$ is nondiagonalizable, there is no inner product in which 
it will be normal.  The defective eigenvalue cannot be on the boundary
of the numerical range~\cite[thm.~1.6.6]{HJ91}.)
If solving a linear system with this $\BA$ via GMRES, $W(\BA)$ will not give
a convergent bound, whereas $W_\BG(\BA)$ will.

\subsection{A preconditioned KKT system}

\begin{figure}[b!]
\begin{center}
\includegraphics[height=2.2in]{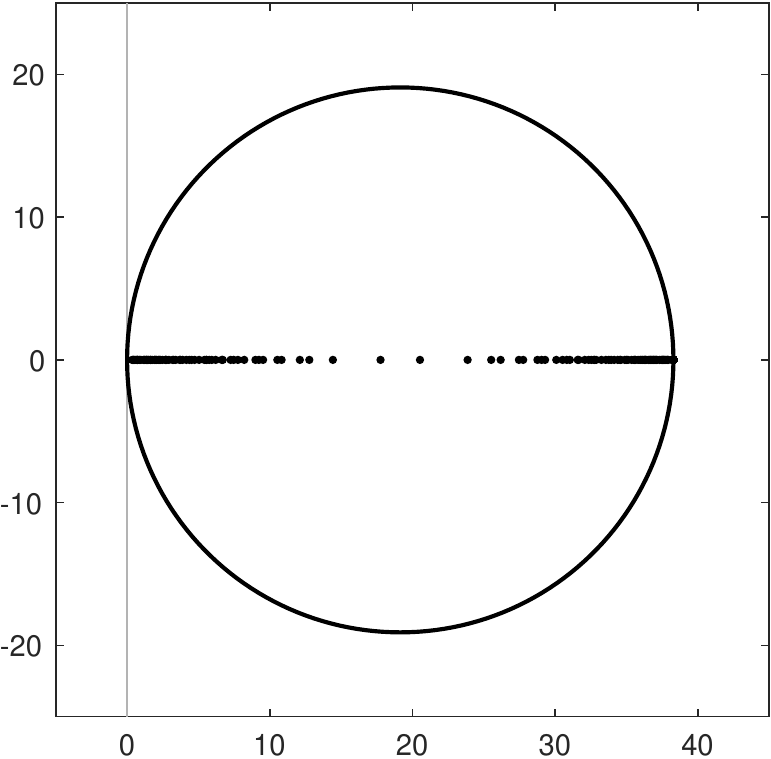}\quad
\includegraphics[height=2.2in]{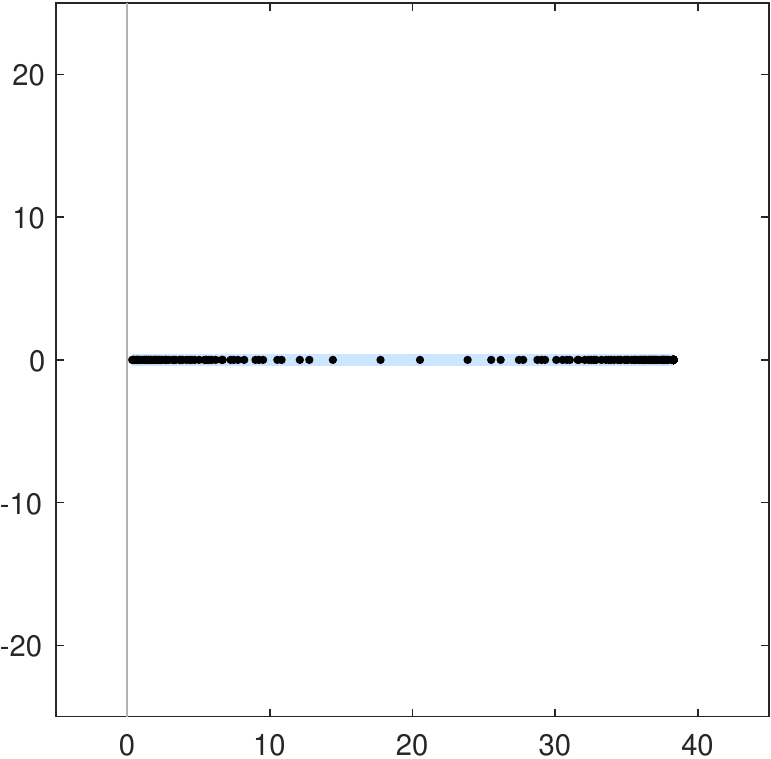}
\end{center}

\begin{picture}(0,0)
\put(118,123){\footnotesize $W(\BA)$}
\put(255,112){\footnotesize $W_\BG(\BA)$}
\end{picture}

\vspace*{-25pt}
\caption{\label{fig:frsw98} 
A preconditioned KKT example, adapted from~\cite{FRSW98}.
Here $\BB\in\R^{64\times 128}$ has normally distributed random entries and 
$\eta = 2@\|\BB\|_2+0.1$.  On the left, the numerical range is a large region
touching the origin: $W(\BA)\cap\R = [0,\eta]$.
On the right, $\BW_\BG(\BA)$ is simply the convex hull of the spectrum, 
$W_\BG(\BA) \approx [0.35633,\eta]$, with $\smash{\sqrt{\kappa_2(\BG)}} \approx 27.6528$.}
\end{figure}

Fischer, Ramage, Silvester, and Wathen~\cite{FRSW98} consider a KKT system with 
the coefficient matrix preconditioned into the form
\[ \BA = \left[\begin{array}{cc}
                \eta@\BI & \BB^* \\[3pt] -\BB & \Bzero  
         \end{array}\right],\]
where $\eta>0$ is a tunable parameter and $\BB\in\R^{M\times N}$.
In the Euclidean inner product the boundary of $W(\BA)$ contains the origin,
since $\BA$ is real-valued and 
\[ \textstyle{\frac{1}{2}} \big(\BA + \BA^*\big)
     = \left[\begin{array}{cc} \eta@\BI & \Bzero \\ \Bzero & \Bzero\end{array}\right].\]
Thus Elman's bound does not apply to such $\BA$; 
see the left plot in \Cref{fig:frsw98} for an example.
Consider the case of $M\le N$ with ${\rm rank}(\BB) = M$.
Fischer et al.~\cite{FRSW98} determine that $\BA$ has a real spectrum
provided $\eta \ge 2@\|\BB\|_2$, and that if $\eta > 2@\|\BB\|_2$, then 
$\BA$ is normal in the inner product defined by 
\[ \BG = \left[\begin{array}{cc} \BI & (2/\eta)@\BB^* \\[3pt] 
               (2/\eta)@\BB & \BI \end{array}\right].
\]
We briefly interpret these results in terms of our Lyapunov setting.
This $\BG$ is positive definite provided $\eta > 2@\|\BB\|_2$, in which
case it defines an inner product in which $\BA$ is self-adjoint
(as $\BA^*\BG = \BG\BA$).  For this $\BG$,
\[ \BA^*\BG+\BG\BA = \left[\begin{array}{cc} 
      2@\eta@\BI - (4/\eta)@\BB^*\BB & 2@\BB^* \\[3pt]
       2@\BB & (4/\eta)\BB\BB^*
     \end{array}\right] =: \BC,
\]
which will be positive definite when $\eta> 2@\|\BB\|_2$
(as this condition ensures the existence of a Cholesky factorization for $\BC$;
alternatively, this condition ensures that $\BA$ is self-adjoint in the $\BG$ inner product 
and all its eigenvalues are positive~\cite[Lemma~2.2]{FRSW98}, so 
$W_\BG(\BA)$ is a real positive interval and $\BC$ must be positive definite).
See the right plot in \Cref{fig:frsw98} for an example.
Indeed, in this case one could run a short-recurrence Krylov method such as 
conjugate gradients or MINRES in the $\BG$ inner product~\cite{FRSW98},
instead of GMRES.

\section{Choosing $\BC$ via Lyapunov inverse iteration}

Different choices of the Hermitian positive definite right-hand side $\BC$
in~\eqref{eq:lyap} give different inner products, and hence different GMRES bounds.
In the absence of application-specific information as in the last section, 
how might one choose $\BC$ to get good bounds?

We first note that this general Lyapunov framework subsumes several 
familiar GMRES bounds.
(i)~Taking $\BC := \BA+\BAs^*$ gives a Lyapunov equation~\cref{eq:lyap} 
with the trivial solution $\BG=\BI$: in this case, \cref{eq:main}
simply recovers the GMRES expression~\cref{eq:bound0}.
(ii)~Suppose $\BA$ is diagonalizable, $\BA=\BV\BLambda\BV^{-1}$, and define
$\BC := \BV^{-*}(\BLambda^*+\BLambda)\BV^{-1}$.  This $\BC$ is
positive definite, since it is a congruence transformation of the positive
definite matrix $\BLambda^*+\BLambda = 2@@\Re(\BLambda)$ 
(via Sylvester's law of inertia~\cite[thm.~4.5.8]{HJ13}, and the 
assumption that the eigenvalues of $\BA$ are in the right half-plane).
The Lyapunov equation $\BAs^*\BG+\BG\BA = \BC$ has the solution
$\BG = \BV^{-*}\BV^{-1}$, so $\sqrt{\kappa_2(\BG)} = \kappa_2(\BV)$, and 
\Cref{thm:main} recovers the standard diagonalization bound for GMRES~\cite{EES83,Elm82}
 for these $\BA$:
\begin{eqnarray*}
\frac{\|\Br_k\|_2}{\|\Bb\|_2} 
       \,\le\, \kappa_2(\BV) \,\mingmres \normG{@p(\BA)} 
       &=& \kappa_2(\BV) \,\mingmres \|\BV^{-1}@p(\BA)\BV\|_2 \\[3pt]
       &=& \kappa_2(\BV) \, \mingmres \max_{\lambda \in \sigma(\BA)} |p(\lambda)|,
\end{eqnarray*}
where $\sigma(\BA)$ denotes the set of eigenvalues of $\BA$.\ \ 
(For the analogous observation about this $\BC$ for 
bounding $\|\eop^{t\BA}\|$, see~\cite[sec.~4.2.2]{Gar23}.)
Of course, for cases where $\kappa_2(\BV)$ is very large, one hopes that a different 
choice for $\BC$ will yield a smaller constant $\sqrt{\kappa_2(\BG)}$, 
while giving a set $W_\BG(\BA)$ that is not much larger 
than the convex hull of $\sigma(\BA)$ (or at least well-separated
from the origin).

One natural choice is $\BC=\BI$ (as used by 
Veseli\'c~\cite{Ves03a},\cite[sect.~21.3]{Ves11}, to bound $\|\eop^{t\BA}\|$),
which we illustrate for two numerical examples in the next section.
As with any positive definite $\BC$, this choice will yield an inner product
in which $W_\BG(\BA)$ is contained in the right half-plane.
To produce other such inner products, one might then iterate:  
Set $\BG_0 = \BC = \BI$ and solve, for $m=1,2,\ldots$
\begin{equation} \label{eq:lyap_invit}
      \BAs^*\BG_{m} + \BG_{m} \BAs = \BG_{m-1}. 
\end{equation}
Each $\BG_m$ gives a new inner product.  
The following simple code will
run one step of this iteration, assuming one starts with an initial
\verb|R = chol(C)|.\\[4pt]
\hspace*{15pt} \verb|R = lyapchol(-A',R');  % Lyapunov solution in factored form| \\
\hspace*{15pt} \verb|A_Gm = R*(A/R);        % Note: 2-norm of A_Gm = G_m-norm of A| \\[4pt]
For the example shown in the next section, 
as $m$ increases $W_{\BG_m}(\BA)$ tends to shrink while $\kappa(\BG_m)$ grows.
For insight into this iteration, consider 
the Lyapunov operator $\Lyap: \C^{n\times n} \to \C^{n\times n}$, 
defined via $\Lyap \BX := \BAs^*\BX+\BX\BA$.\ \ 
Then the recursion \cref{eq:lyap_invit} amounts
to inverse iteration: $\BG_{m} = \Lyap^{-1} \BG_{m-1}$.  
Standard convergence theory for inverse iteration gives insight into the
behavior of $\BG_m$ as $m\to\infty$.  Suppose $\BA$ is diagonalizable with 
eigenvalues $\{\lambda_j\}_{j=1}^n$ and corresponding \emph{left} eigenvectors $\{\Bw_j\}_{j=1}^n$,
i.e., $\Bw_j^*\BA = \lambda_j\Bw_j^*$.
Then $\Lyap (\Bw_j^{}\Bw_\ell^*) = (\overline{\lambda_j} +\lambda_\ell) (\Bw_j^{}\Bw_\ell^*)$,
and so $\overline{\lambda_j} +\lambda_\ell$ is an eigenvalue of $\Lyap$~\cite[sect~4.4]{HJ91}.
If we assume for simplicity that $\BA$ has distinct, real, positive eigenvalues 
$\lambda_1 < \lambda_2 < \cdots < \lambda_n$, then the smallest magnitude 
eigenvalue of $\Lyap$ will be $2@\lambda_1$ with corresponding eigenvector $\Bw_1^{}\Bw_1^*$.  
We expect $\BG_m$ to converge toward a multiple of this rank-1 matrix, suggesting the
significance of $\lambda_1$ in $W_{\BG_m}(\BA)$ and why $\kappa(\BG_m)$ will increase.
(For other applications of Lyapunov inverse iteration, see~\cite{EMSW12,GMOXZ06,MS10}.)

This inverse iteration perspective suggests another way to generate $\BC$.\ \ 
Since $\BA-s@\BI$ has all eigenvalues in the right half-plane for any real $s>0$ 
that is smaller than the real part of the leftmost eigenvalue of $\BA$, 
one could iterate
\begin{equation} \label{eq:lyap_shifted}
      (\BAs-s@\BI)^*\BG_{m} + \BG_{m} (\BAs-s@\BI) = \BG_{m-1}.
\end{equation}
This approach amounts to Lyapunov shifted inverse iteration, since 
$\Lyap\BX - 2@s\BX = (\BA-s@\BI)^*\BX + \BX(\BA-s@\BI)$;
thus $\BG_{m} = (\Lyap - 2@s@\Id)^{-1}\BG_{m-1}$,
where $\Id$ is the identity map on $\C^{n\times n}$.
With \verb|I = eye(n)|, a step of this iteration can be implemented as follows.\\[4pt]
\hspace*{15pt} \verb|R = lyapchol(-(A-s*I)',R'); % Lyapunov sol'n in factored form| \\
\hspace*{15pt} \verb|A_Gm = R*(A/R);             % Note: norm(A_Gm) = G_m-norm of A| 

\medskip
More generally, one might expect different choices of $\BC$ to give
better bounds on $\|\Br_k\|_2$ for different values of $k$: 
early iterations would benefit from a small value of $\kappa_2(\BG)$, 
even if that meant a slow convergence rate; 
later iterations might benefit from an inner product 
having a larger $\kappa_2(\BG)$ but a smaller $W_\BG(\BA)$.
Finding an \emph{optimal} choice for $\BC$ appears to be quite difficult
(for analogous optimization for bounding $\|\eop^{t\BA}\|$ as a
function of $t$, see~\cite[sect.~4.2]{Gar23});
however, the examples in the next section suggest that the iterations~\cref{eq:lyap_invit} 
and~\cref{eq:lyap_shifted} can generate a range of useful inner products.

\section{Examples of Lyapunov inverse iteration}

\begin{figure}[b!]
\begin{center}
\includegraphics[height=1.7in]{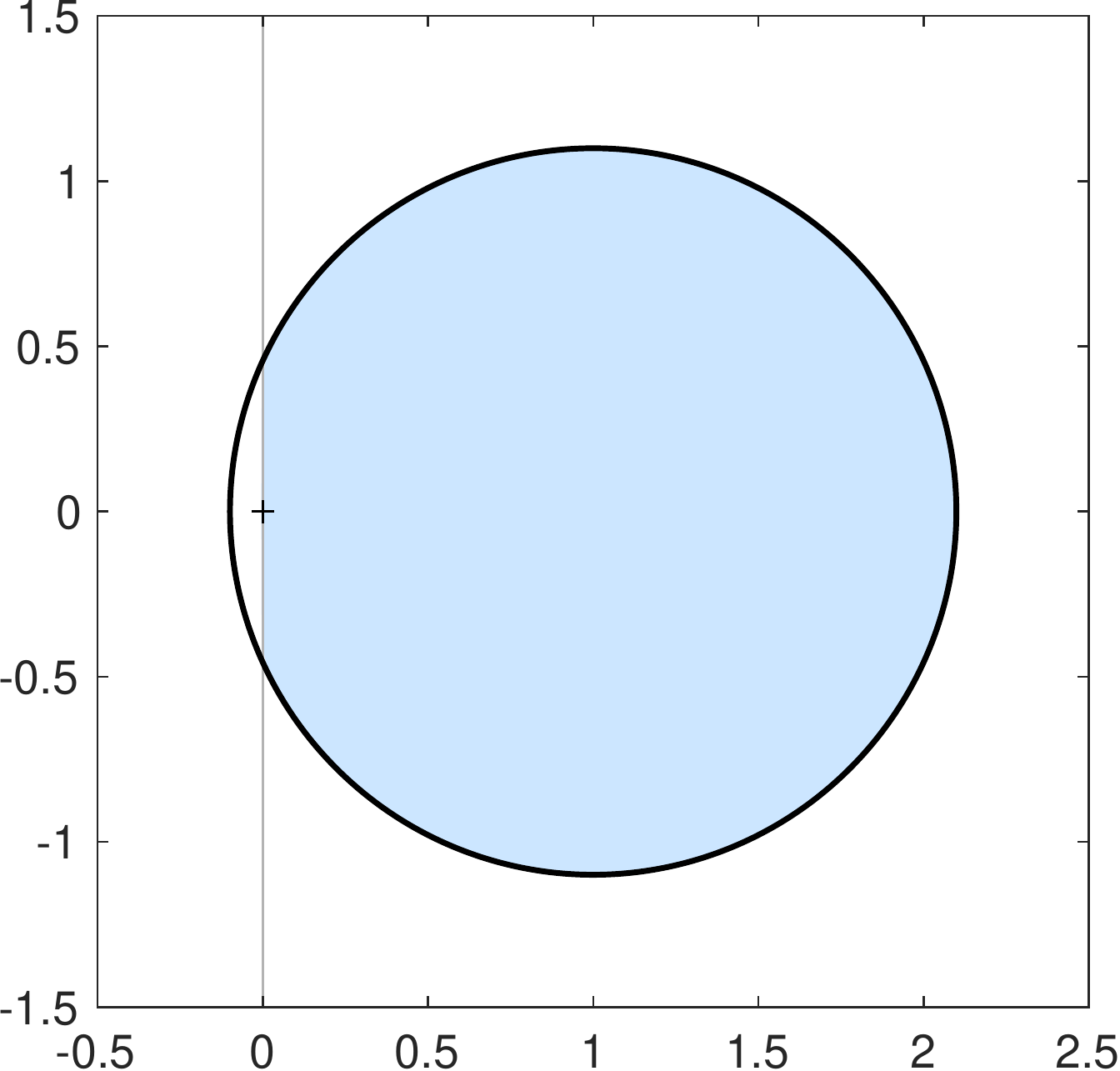}\hspace*{30pt}
\includegraphics[height=1.7in]{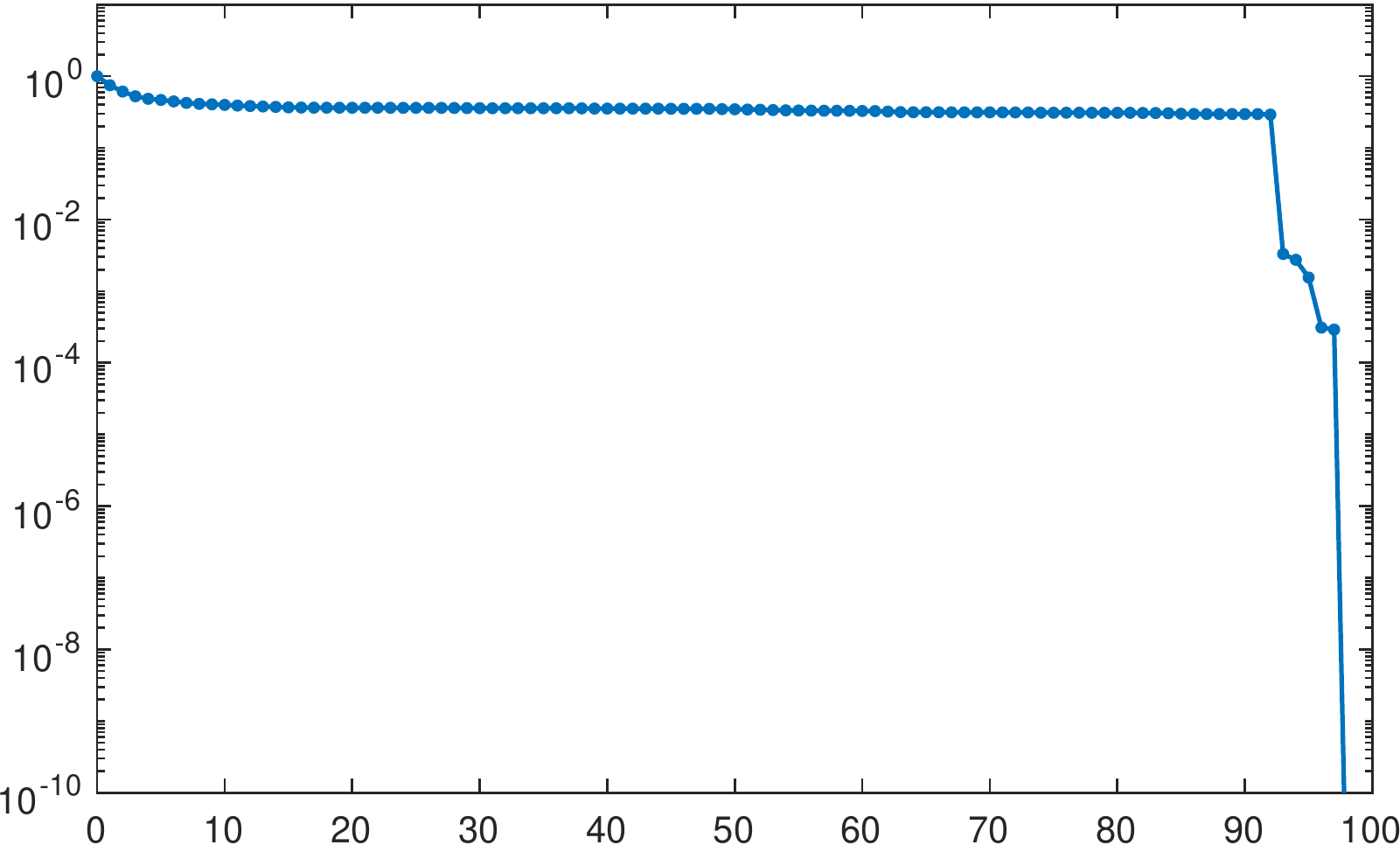}
\end{center}

\begin{picture}(0,0)
\put(58,84){\footnotesize $W_\BG(\BA)$}
\put(140,120){\footnotesize $\displaystyle{\frac{\|\Br_k\|_2}{\|\Bb\|_2}}$}
\put(238,10){\footnotesize iteration, $k$}
\put(23,44){\vector(1,4){7.5}}
\put(17,33){\whiteout{20pt}{10pt}}
\put(17,36){\footnotesize $W(\BA)$}
\end{picture}

\vspace*{-14pt}
\caption{\label{fig:jordan_ex}
Jordan block example.  
Solving $\BAs^*\BG+\BG\BA = \BI$ yields an inner product
in which $W_\BG(\BA)$ (blue region) is in the right half-plane, but barely:
$\mu_\BG(\BA)\approx 3.6\times 10^{-9}$, 
with $\smash{\sqrt{\kappa(\BG)}} \approx 2.4\times 10^4$.
These quantities are consistent with the slow convergence of GMRES observed on 
the right (for 100 random $\Bb$ vectors).}
\end{figure}

Two simple examples illustrate how this Lyapunov inner product approach 
can yield insight about GMRES convergence.  
Consider the two bidiagonal matrices of dimension $n=100$ given by
\[ \BA = \left[\begin{array}{ccccc} 1 & \alpha \\ & 1 & \alpha \\ 
          & & 1 & \ddots \\ & & & \ddots & \alpha \\ & & & & 1 \end{array}\right], \qquad
   \BA = \left[\begin{array}{ccccc} 1 & \gamma/1 \\ & 1 & \gamma/2 \\ 
          & & 1 & \ddots \\ & & & \ddots & \gamma/(n-1) \\ & & & & 1 \end{array}\right].
\]
With $\alpha=1.1$ and $\gamma=2$ the numerical range contains the origin for 
both matrices.
The matrix on the left is a scaled Jordan block with eigenvalue $\lambda=1$; 
the matrix on the right also has $\lambda=1$ as a maximally defective eigenvalue,
but the decreasing superdiagonal entries lead to interesting GMRES convergence,
as noted in~\cite{DTT98} (where $\BA$ is called an ``integration matrix'').
The pseudospectra of these two matrices differ considerably, and can give 
insight into GMRES convergence~\cite{DTT98,Emb22}.
In both cases, we solve $\BAs^*\BG+\BG\BA = \BI$ to get
an inner product in which $W_\BG(\BA)$ is contained in the right half-plane.
For the Jordan block, \Cref{fig:jordan_ex} compares $W(\BA)$ (bounded by a black line) 
to $W_\BG(\BA)$ (blue region), and shows GMRES convergence for 
100 different $\Bb$ vectors having normally distributed random entries.  
The Lyapunov solution gives an inner product in which $W_\BG(\BA)$
is in the right half-plane, but only barely; $\sqrt{\kappa(\BG)}$ is very large.  
Both of these properties are consistent with the slow convergence of GMRES.

\begin{figure}[t!]
\begin{center}
\includegraphics[height=1.7in]{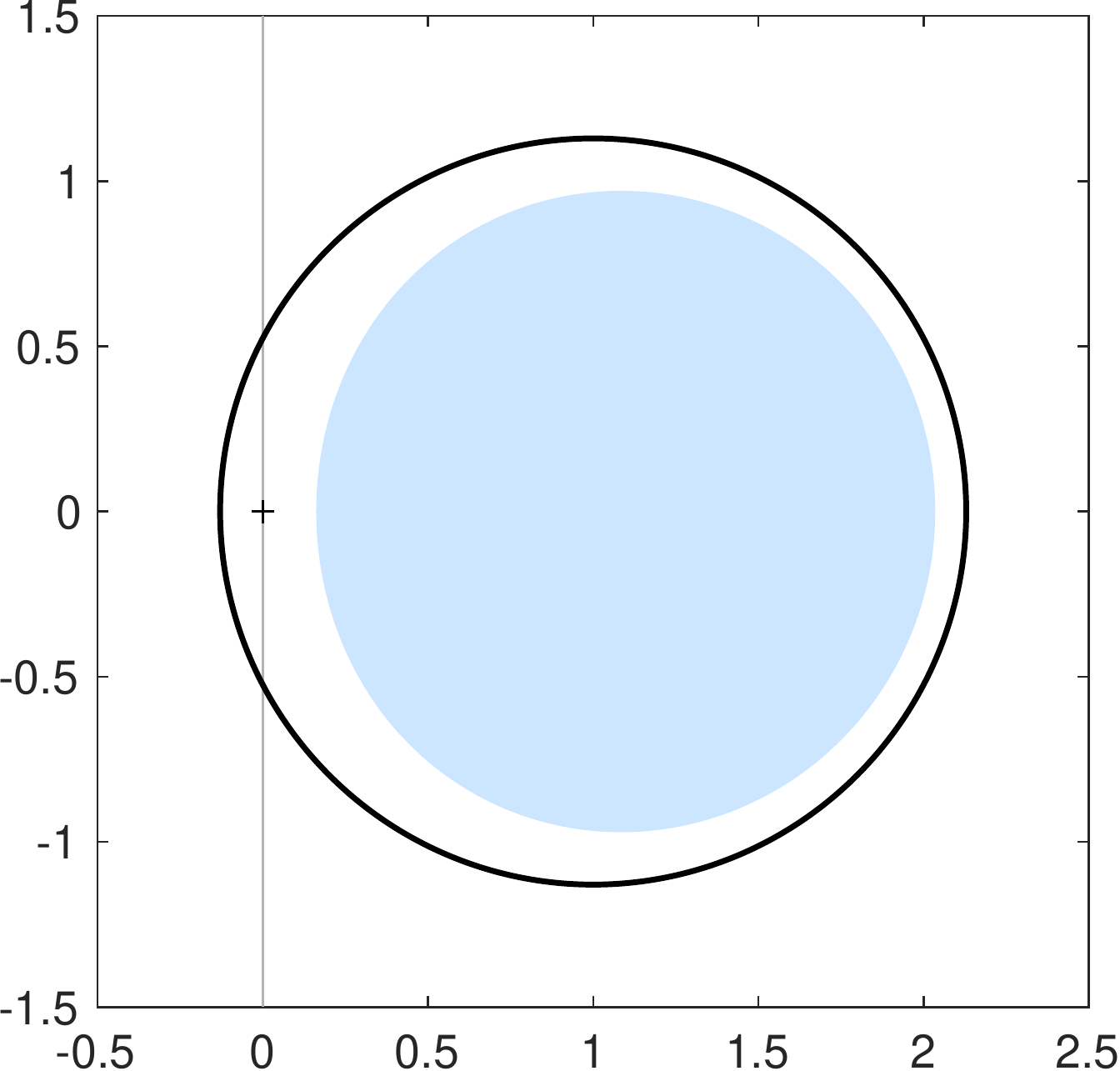}\hspace*{30pt}
\includegraphics[height=1.7in]{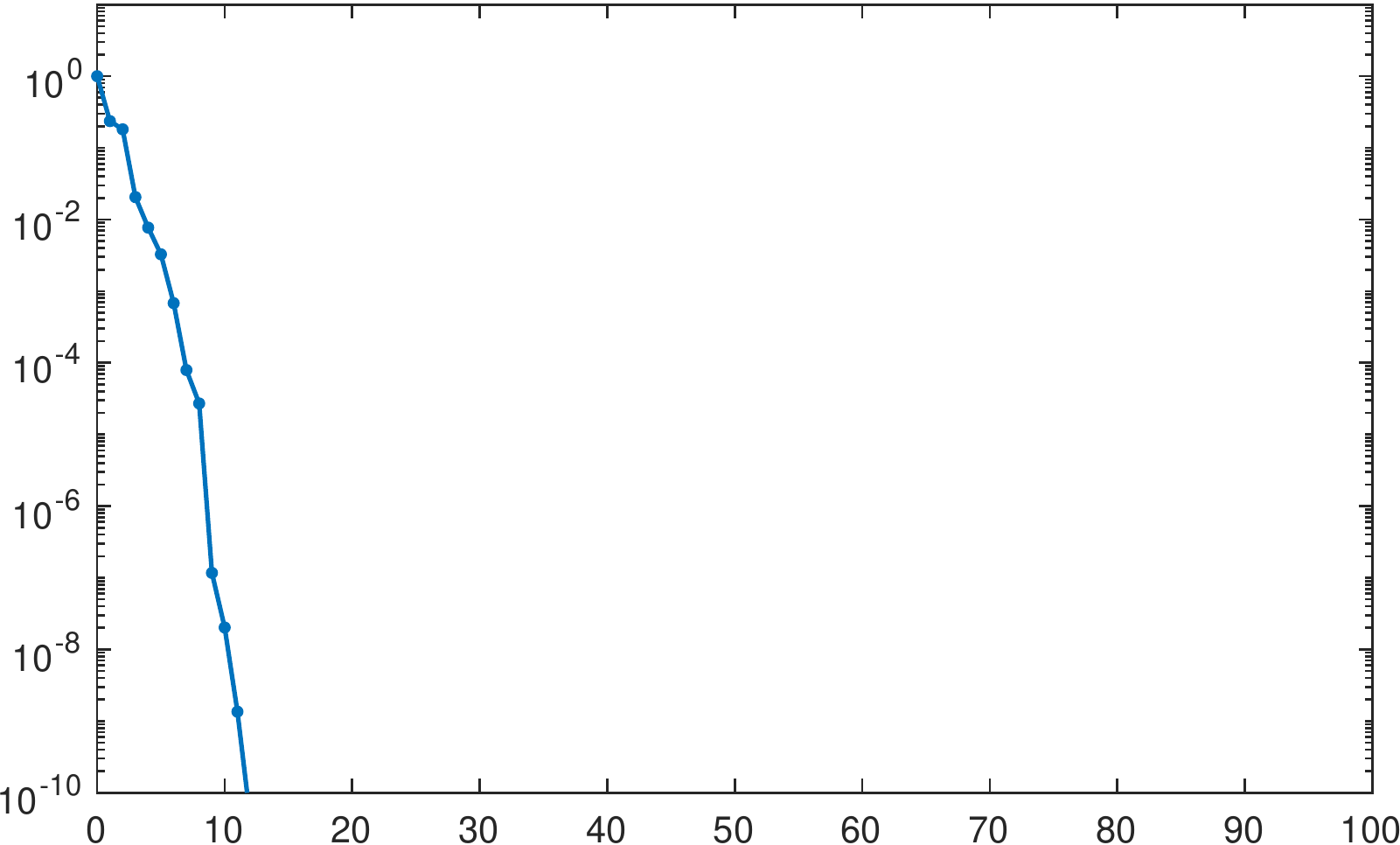}
\end{center}

\begin{picture}(0,0)
\put(58,84){\footnotesize $W_\BG(\BA)$}
\put(140,120){\footnotesize $\displaystyle{\frac{\|\Br_k\|_2}{\|\Bb\|_2}}$}
\put(238,10){\footnotesize iteration, $k$}
\put(23,44){\vector(1,4){7}}
\put(17,33){\whiteout{20pt}{10pt}}
\put(17,36){\footnotesize $W(\BA)$}
\end{picture}

\vspace*{-14pt}
\caption{\label{fig:intmat_ex}
Integration matrix example.  
Solving $\BAs^*\BG+\BG\BA = \BI$ yields an inner product
in which $W_\BG(\BA)$ (blue region) is well separated from the origin:
$\mu_\BG(\BA)\approx 0.17$, with $\smash{\sqrt{\kappa(\BG)}} \approx 3.5$.
These quantities are consistent with the rapid convergence of 
GMRES observed on the right (for 100 random $\Bb$ vectors).}
\end{figure}

\begin{figure}[t!]
\hspace*{18pt}\includegraphics[width=1.8in]{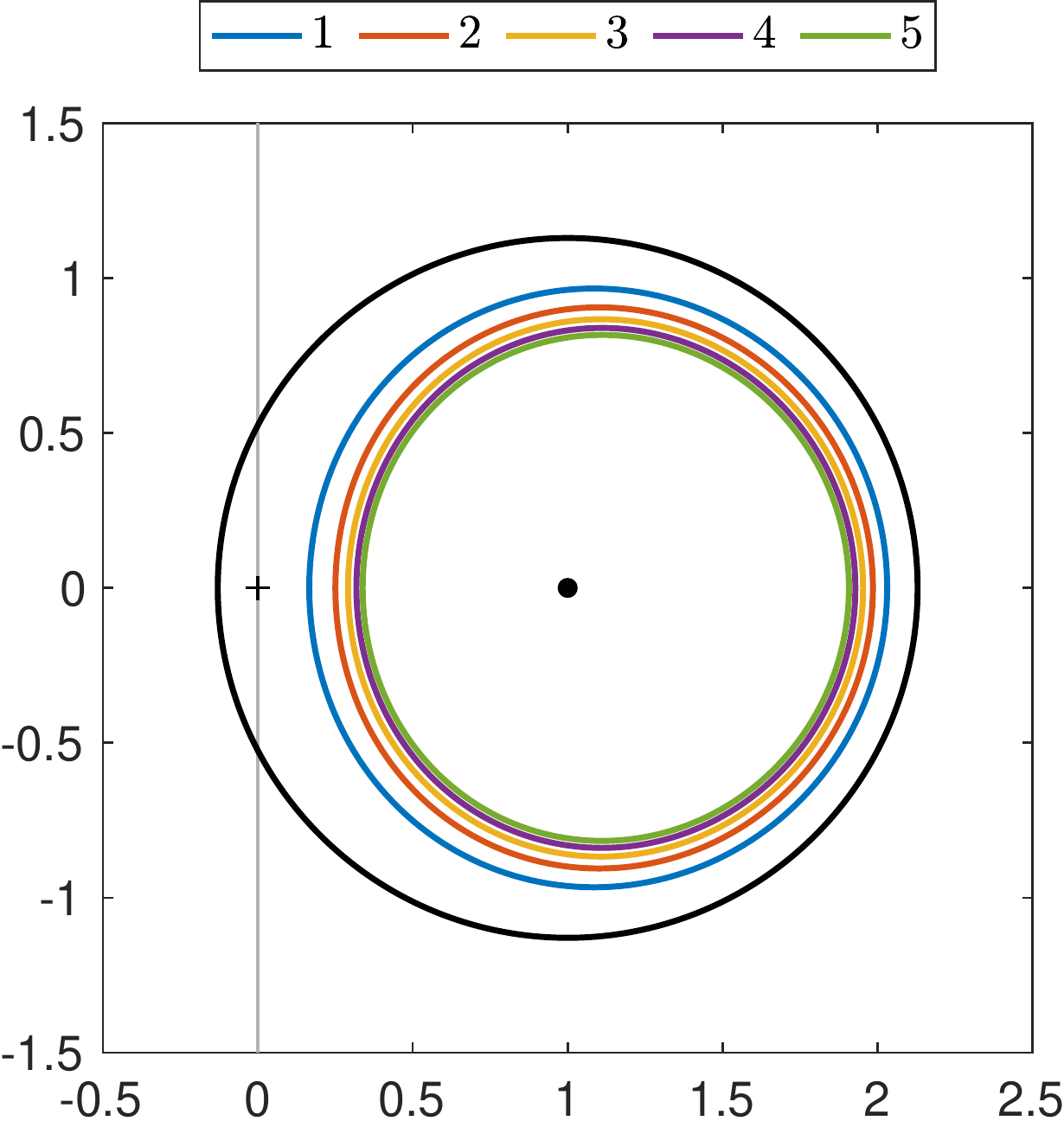}
\begin{picture}(0,0)
 \put(-120,131){\scriptsize $m$:}
\end{picture}
\hspace*{20pt}
\raisebox{62pt}{{\small \begin{tabular}{crcc}
$m$ &   \multicolumn{1}{c}{$\sqrt{\kappa(\BG_m)}$} & $\mu_{\BG_m}(\BA)$ & $\|\BA\|_{\BG_m}$ \\ \hline
 1 &     3.49787 &      0.16600  &     2.21253\\ 
 2 &     9.21667 &      0.25027  &     2.12643\\ 
 3 &    21.34399 &      0.29110  &     2.07321\\ 
 4 &    45.58853 &      0.31813  &     2.03461\\ 
 5 &    91.87710 &      0.33835  &     2.00391\\ 
\end{tabular}}}

\vspace*{-3pt}
\caption{\label{fig:iterate}
Using Lyapunov inverse iteration to construct inner products:
$\BAs^*\BG_{m} + \BG_{m}\BA = \BG_{m-1}$ with $\BG_0=\BI$,
applied to the integration matrix.
The plot shows the boundary of $W_{\BG_k}(\BA)$.
(The outermost (black) line indicates $W(\BA)$; 
the $m=1$ line matches the blue region in \Cref{fig:intmat_ex}.)}
\end{figure}

\begin{figure}[t!]
\hspace*{18pt}\includegraphics[width=1.8in]{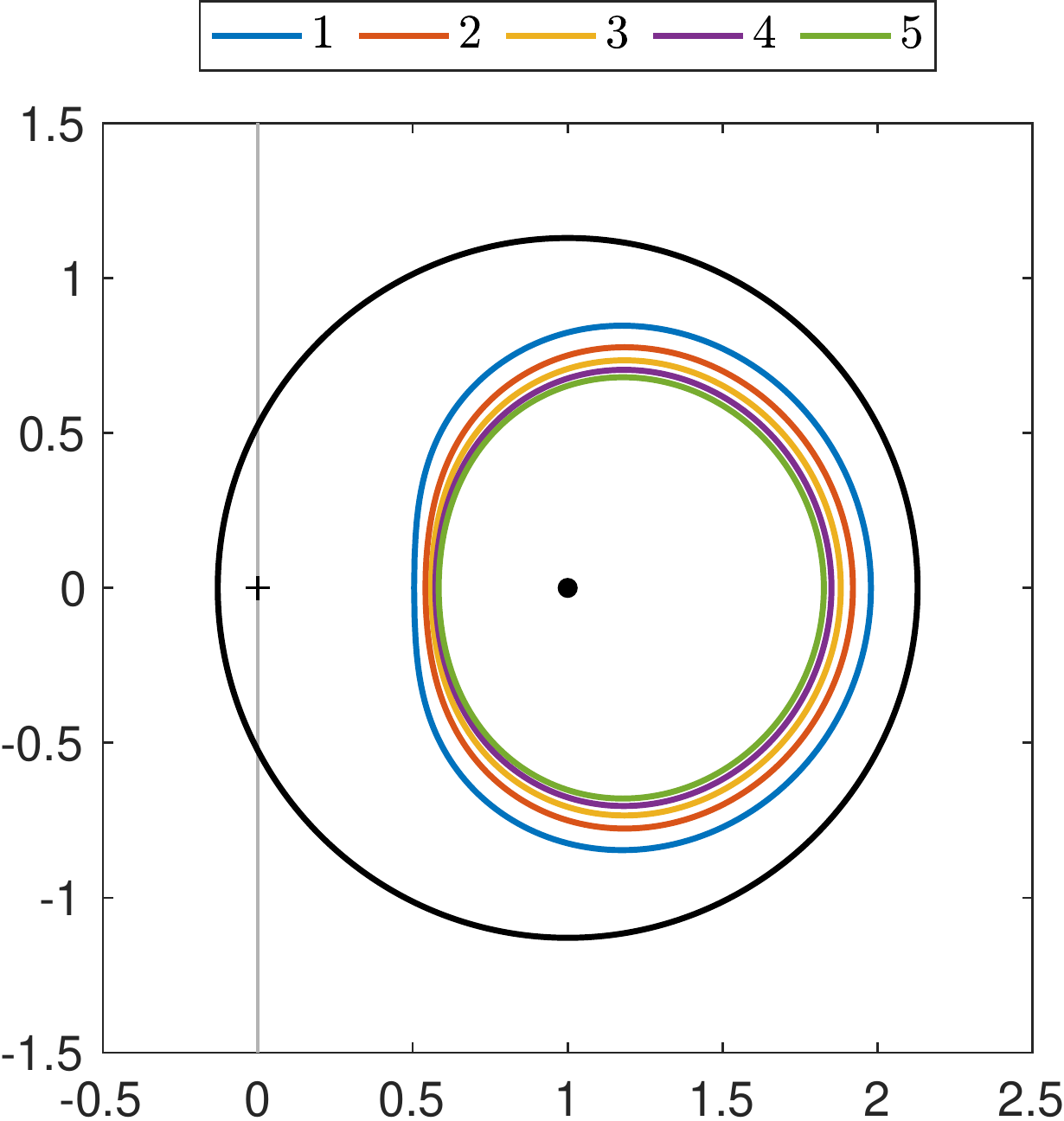}
\begin{picture}(0,0)
 \put(-120,131){\scriptsize $m$:}
\end{picture}
\hspace*{20pt}
\raisebox{62pt}{{\small \begin{tabular}{crcc}
$m$ &   \multicolumn{1}{c}{$\sqrt{\kappa(\BG_m)}$} & $\mu_{\BG_m}(\BA)$ & $\|\BA\|_{\BG_m}$ \\ \hline
 1 &    18.44026 &      0.50435  &     2.11270\\ 
 2 &    86.38039 &      0.54103  &     2.02089\\ 
 3 &   312.50604 &      0.55974  &     1.96449\\ 
 4 &   980.52145 &      0.57312  &     1.92343\\ 
 5 &  2791.92538 &      0.58366  &     1.89110\\ 
\end{tabular}}}

\vspace*{-3pt}
\caption{\label{fig:shift}
Using shifted Lyapunov inverse iteration to construct inner products:
$(\BAs^*-s@\BI) \BG_{m} + \BG_{m}(\BA-s@\BI) = \BG_{m-1}$ with $\BG_0=\BI$,
applied to the integration matrix with shift $s=0.5$.
The plot shows the boundary of $W_{\BG_m}(\BA)$.
(The outermost (black) line indicates $W(\BA)$.)}
\end{figure}

In contrast, GMRES converges much faster for the integration matrix,
as seen in \Cref{fig:intmat_ex}.
Now $W_\BG(\BA)$ is much better separated from the origin,
and $\sqrt{\kappa(\BG)}$ is quite modest, 
consistent with the rapid convergence of GMRES.\ \   
\Cref{fig:iterate} shows $W_{\BG_m}(\BA)$ for inner products 
generated by Lyapunov inverse iteration~\cref{eq:lyap_invit} with $\BG_0=\BI$:
the numerical range shrinks while the constant $\sqrt{\kappa(\BG_m)}$ grows.
This effect is even more pronounced for shifted Lyapunov inverse 
iteration~\cref{eq:lyap_shifted}, as seen in \Cref{fig:shift}.

Consider the convergence rates exhibited by GMRES bounds
that are based on these $\BG_m$ inner products: 
the rate $\rho_{\rm E} = (1-\mu_\BG(\BA)^2/\|\BA\|_\BG^2)^{1/2}$ from \Cref{cor:elman}, 
the rate $\rho_\beta$ from~\Cref{cor:bgt}, 
and the rate $\rho_\BG = \max_{z\in W_\BG(\BA)} |1-z/c|$, that comes from bounding $W_\BG(\BA)$
by a circle centered at $c\in\R$, the average of the extreme real values of $W_\BG(\BA)$.
Note that $\rho_\BG$ can be used to get an upper bound for \Cref{cor:cp}, 
since 
\[ \mingmres \max_{z\in W_\BG(\BA)} |p(z)| 
   \le \max_{z\in W_\BG(\BA)} |1-z/c|^k = \rho_\BG^k.\]
\Cref{tbl:rates} shows the rates for several of the numerical ranges
from \Cref{fig:iterate,fig:shift}.  
Even crudely bounding $W_\BG(\BA)$ by a circle gives
significantly improved rates over $\rho_{\rm E}$ and $\rho_\beta$, in this case.

\begin{figure}[h!]
\smallskip
\begin{center}
\includegraphics[height=1.7in]{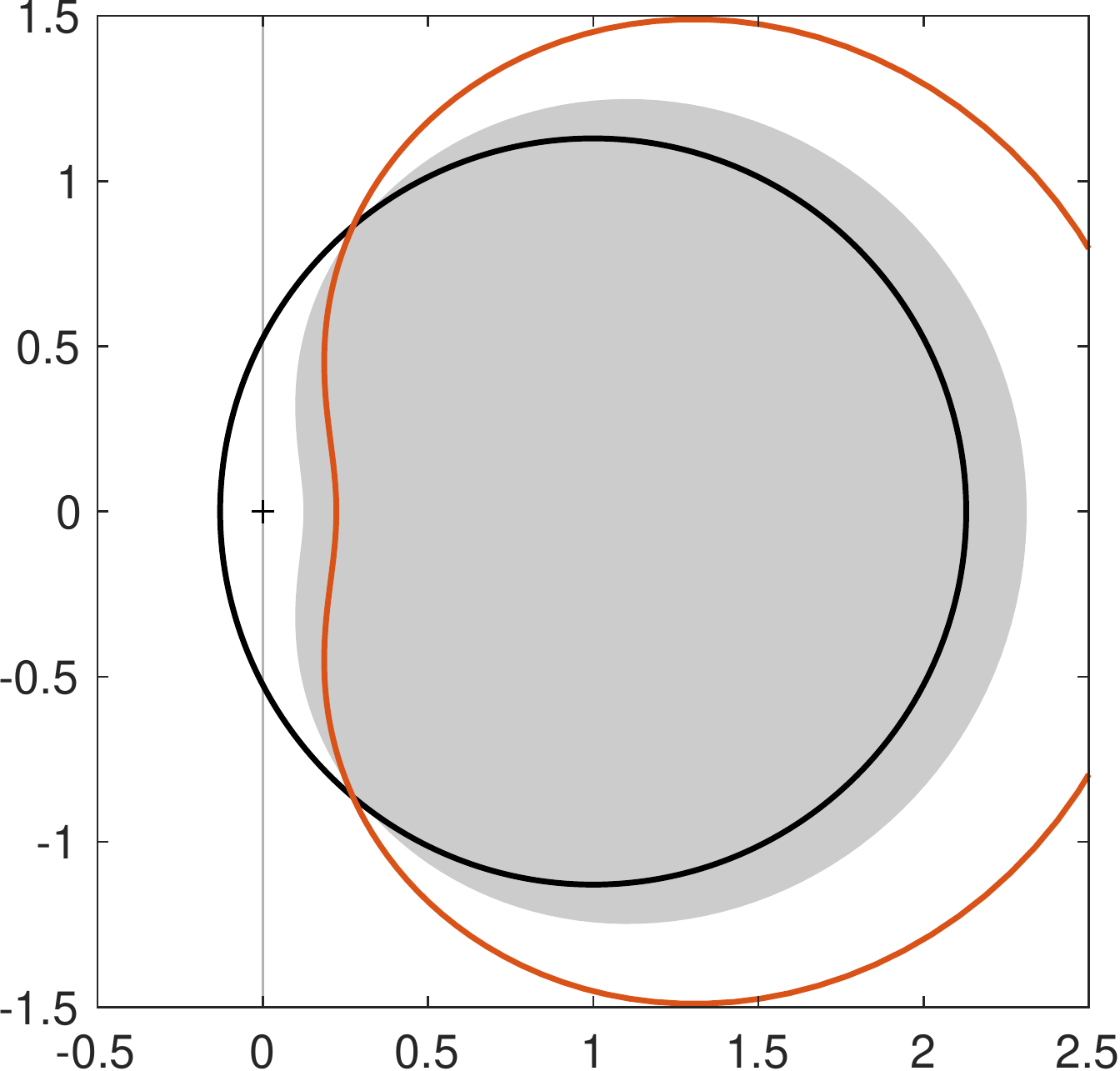}\hspace*{30pt}
\includegraphics[height=1.7in]{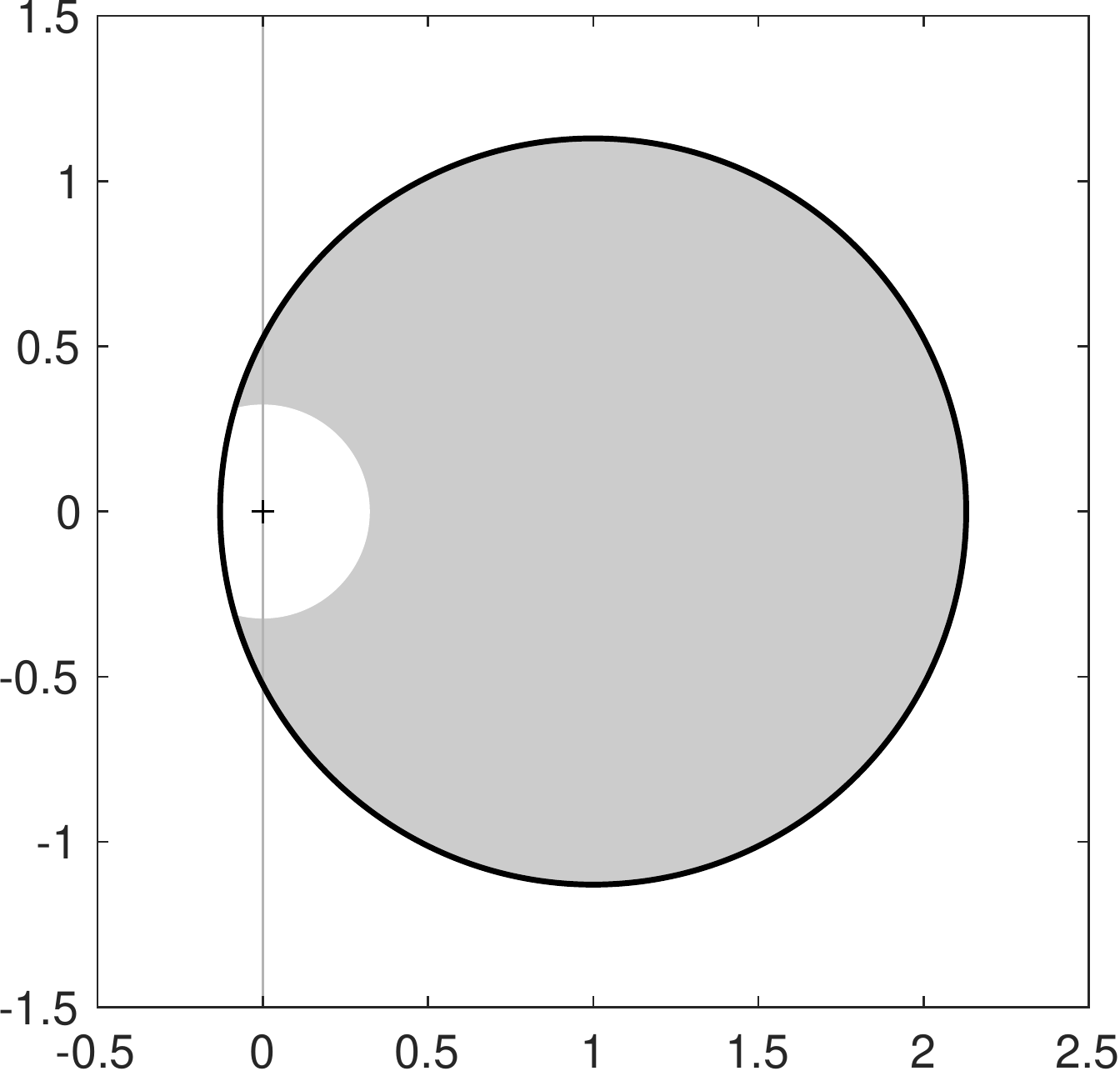}
\end{center}

\begin{picture}(0,0)
\put(83,84){\footnotesize $W(\BAs^{1/2})^2$}
\put(57,38){\vector(1,4){5}}
\put(47,30){\whiteout{20pt}{10pt}}
\put(47,33){\footnotesize $W(\BA)$}
\put(230,70){\footnotesize \shortstack{$W(\BA)$\\$\cap$\\$\{|z|>1/r(\BAs^{-1})\}$}}
\put(219,38){\vector(1,4){5}}
\put(209,30){\whiteout{20pt}{10pt}}
\put(209,33){\footnotesize $W(\BA)$}
\end{picture}

\vspace*{-22pt}
\caption{\label{fig:intmat_ex_CG}
Integration matrix example.  The gray regions show alternatives
to the numerical range that exclude the origin.  
The solid black line shows the boundary of $W(\BA)$, for comparison.
The red line in the left plot indicates the boundary of 
$\exp(W(\log(\BA)))$, another set studied in~\cite{CG15}.}
\end{figure}

Let us compare the numerical ranges for the integration matrix 
shown in \Cref{fig:intmat_ex,fig:iterate,fig:shift} 
to several other sets that excise the origin from the numerical range of $\BA$, 
for purposes of analyzing GMRES.\ \ 
Choi and Greenbaum~\cite{CG15} observe that for nonsingular $\BA$,
 $W(\BAs^{1/d})$ will exclude the origin for $d$ sufficiently large, 
in which case one can apply the Crouzeix--Palencia theorem to obtain
\[ \frac{\|\Br_k\|_2}{\|\Bb\|_2} 
     \,\le \,(1+\sqrt{2}) \,\mingmres \max_{z\in W(\BAs^{1/d})^d} |p(z)|,\]
where $W(\BAs^{1/d})^d := \{ z^d : z \in W(\BAs^{1/d})\}$.
For the integration matrix example, $d=2$ suffices.  \Cref{fig:intmat_ex_CG} (left)
shows that the set $W(\BAs^{1/2})^2$ is in the right half-plane, though it contains
some points outside $W(\BA)$.  
The same is true of the set $\exp(W(\log(\BA))) = \lim_{d\to\infty} W(\BAs^{1/d})^d$ 
also studied in~\cite{CG15}.

\begin{table}[t]
\caption{\label{tbl:rates}
Convergence rates associated with several $W_{\BG_m}(\BA)$ for the integration matrix.
Inverse iteration uses~\cref{eq:lyap_invit}, 
shifted inverse iteration uses~\cref{eq:lyap_shifted} with $s=0.5$; 
both start with $\BG_0=\BI$.}
\begin{center} \small
\begin{tabular}{l|cc|cc}
& \multicolumn{2}{c}{\emph{inverse iteration}} 
& \multicolumn{2}{c}{\emph{shifted inverse iteration}} \\
& \ $m=1$ 
& \ $m=5$ 
& \ $m=1$ 
& \ $m=5$  \\ \hline
$\rho_{\rm E}$  &\ \  0.99718\ \  &\ \ 0.98564\ \  & \ \ 0.97109\ \ &\ \ 0.95118\ \ \\
$\rho_\beta$    &\ \  0.94257\ \  &\ \ 0.87140\ \  & \ \ 0.81859\ \ &\ \ 0.76566\ \ \\
$\rho_\BG$      &\ \  0.88107\ \  &\ \ 0.72816\ \  & \ \ 0.69335\ \ &\ \ 0.56739\ \ \\
\end{tabular}
\end{center}
\end{table}

Crouzeix and Greenbaum~\cite{CG19} describe another technique 
that removes the origin from $W(\BA)$.  Let $r(\BX) = \max_{z\in W(\BX)} |z|$ 
denote the \emph{numerical radius} of the matrix $\BX\in\Cnn$.  The set
\[ \Omega\ :=\ W(\BA)\ \cap\ \{z\in \C: |z| \ge 1/r(\BAs^{-1})\} \]
excludes the origin, 
and the GMRES residual then satisfies~\cite[sect.~7.1]{CG19}
\[ \frac{\|\Br_k\|_2}{\|\Bb\|_2} 
     \,\le \,(3+2\sqrt{3}) \,\mingmres \max_{z\in \Omega} |p(z)|.\]
\Cref{fig:intmat_ex_CG} (right; repeated from~\cite[fig.~3.11]{Emb22}) 
shows the set $\Omega$, which cuts a portion of a circular disk out of $W(\BA)$.

\section{Conclusions}

When all the eigenvalues of the matrix $\BA$ are contained 
in the right half-plane, 
one can solve the Lyapunov equation $\BAs^*\BG + \BG\BA = \BC$ 
with any Hermitian positive definite $\BC$
to obtain an inner product in which the numerical range of $\BA$ 
is also contained in the right half-plane.
One can then use this new inner product to 
analyze the convergence of GMRES 
(run in the standard Euclidean inner product).
Different choices of $\BC$ give
different GMRES bounds, which could better predict convergence
for early or later iterations.  This observation opens several
avenues for exploration, not only for the study of GMRES convergence,
but potentially for the development and analysis of inner products
(preconditioners) in which GMRES might converge more favorably.


\section*{Acknowledgments}

I thank Anne Greenbaum, Pierre Marchand, Volker Mehr\-mann, and Andy Wathen 
for their comments and for suggesting several helpful references.

\end{document}